\documentclass{article}
\usepackage{amssymb}
\usepackage{amsmath}
\usepackage[hyperindex,breaklinks,colorlinks,citecolor=blue,pagebackref]{hyperref}
\usepackage[pdftex]{graphicx}

\newtheorem{theorem}{Theorem}
\newtheorem{acknowledgement}[theorem]{Acknowledgement}

\newtheorem{corollary}[theorem]{Corollary}

\newtheorem{definition}[theorem]{Definition}

\newtheorem{lemma}[theorem]{Lemma}

\newtheorem{problem}[theorem]{Problem}
\newtheorem{proposition}[theorem]{Proposition}
\newtheorem{remark}[theorem]{Remark}

\newenvironment{proof}[1][Proof]{\noindent\textbf{#1.} }{\ \rule{0.5em}{0.5em}}
\input{tcilatex}
\begin{document}

\title{Dynamical systems, simulation, abstract computation.}
\author{Stefano Galatolo \and Mathieu Hoyrup \and Crist\'obal Rojas}
\maketitle

\begin{abstract}
We survey an area of recent development, relating dynamics to theoretical
computer science. We discuss the theoretical limits of simulation and
computation of interesting quantities in dynamical systems. We will focus on
central objects of the theory of dynamics, as invariant measures and
invariant sets, showing that even if they can be computed with arbitrary
precision in many interesting cases, there exists some cases in which they
can not. We also explain how it is possible to compute the speed of
convergence of ergodic averages (when the system is known exactly) and how
this entails the computation of arbitrarily good approximations of points of
the space having typical statistical behaviour (a sort of constructive
version of the pointwise ergodic theorem).
\end{abstract}

\tableofcontents

\section{Introduction}

The advent of automated computation led to a series of great successes and
achievements in the study of dynamical problems.

The use of computers and simulations allowed to compute and forecast the
behavior of many important natural phenomena and, on the other hand, led to
the discovery of important general aspects of dynamics.

This motivates the huge work that was made by hundred of scientists to
improve \textquotedblleft practical\textquotedblright\ simulation and
computation techniques.

It also motivates the study of the theoretical limits of simulation and
computation techniques, and the theoretical understanding of related
problems.

In this paper we want to focus on some of these aspects related to \emph{%
rigorous }computation and simulation of (discrete time) dynamical systems.

The simulation and investigation of dynamics started with what we call the
``naive'' approach, in which the user just implements the dynamics without
taking rigorous account of numerical errors and roundings. Then he ``looks''
to the screen to see what happens.

Of course, the sensitivity to initial conditions, and the typical
instability of many interesting systems (to perturbations on the map
generating the dynamics) implies that what it is seen on the screen could be
completely unrelated to what was meant to be simulated.

In spite of this, the naive approach turns out to work \textquotedblleft
unreasonably\textquotedblright\ well in many situations and it is still the
most used one in simulations. The theoretical reasons why this method works
and its limits, in our opinion are still to be understood (some aspects have
been investigated in \cite{Lax},\cite{mie},\cite{Gr} \cite{Bla94, Bla89}
e.g.).

On the opposite side from the naive approach, there is the \textquotedblleft
absolutely rigorous\textquotedblright\ approach, which will \ be the main
theme of this paper: the user looks for an algorithm which can give a
description of the object which is meant to be computed, up to any desired
precision.

In this point of view, for example the constant $e$ is a computable number
because there is an algorithm that is able to produce a rational
approximation of $e$ at any given precision (for example finding the right $%
m $ and calculating $\sum_{1}^{m}\frac{1}{n!}$ such that the error is
smaller than requested).

In an \textquotedblleft absolutely rigorous\textquotedblright\ simulation
(computation) the initial point or some initial distribution is supposed to
be known up to any approximation and the transition map is also supposed to
be known up to any accuracy (suitable precise definitions will be given
below). This allows the evolution of the system to be simulated with any
given accuracy, and the question arise, if interesting and important objects
related to the dynamics (invariant sets or invariant measures e.g.) can be
computed from the description of the system (or maybe adding some additional
information).

In this paper we would like to give a survey of a group of results related
to these computational aspects, restating and updating them with some new
information. We will see that in many cases the interesting objects can be
computed, but there are some subtleties, and cases where the interesting
object cannot be computed from the description of the system or cannot be
computed at all (again, up to any given precision). In particular, this
happen case for the computation of invariant measures and Julia sets.

Hence, these results set theoretical limits to such computations.

\medskip \noindent \textbf{Computing invariant measures.}

An important fact motivating the study of the statistical properties of
dynamical systems is that the pointwise long time prediction of a chaotic
system is not possible, whereas, in many cases, the estimation or
forecasting of averages and other long time statistical properties is. This
often corresponds in mathematical terms to computing invariant measures, or
estimating some of their properties, as measures contain information on the
statistical behavior of the system $(X,T)$ and on the potential behavior of
averages of observables along typical trajectories of the system (see
Section \ref{DS}).

An invariant measure is a Borel probability measure $\mu $ on $X$ such that
for each measurable set $A$ it holds $\mu (A)=\mu (T^{-1}(A))$. They
represent \emph{equilibrium} states, in the sense that probabilities of
events do not change in time.

Rigorously, compute an invariant measure means to find an algorithm which is
able to output a description of the measure (for example an approximation of
the measure made by a combination of delta measures placed on
\textquotedblleft rational\textquotedblright\ points) up to any prescribed
precision.

We remark that once an interesting invariant measure is computed, it is
possible to deduce from it several other important information about the
dynamics: Lyapunov exponents, entropy, escape rates, etc... For example, in
dimension one, once an ergodic invariant measure $\mu $ has been
approximated, the Lyapunov exponent $\lambda _{\mu }$ can be estimated using
the formula $\lambda _{\mu }=\int_{0}^{1}\log _{2}T^{\prime }d\mu $, where $%
T^{\prime }$ denotes the derivative of the map generating the dynamics. In
higher dimensions, similar techniques can be applied (see e.g. \cite{F08}
for more examples of derivation of dynamical quantities from the computation
of the invariant measure).

Before giving more details about the computation of invariant measures, we
remark that, since there are only countably many \textquotedblleft
algorithms\textquotedblright\ (computer programs), whatever we mean by
\textquotedblleft approximating a measure by an algorithm\textquotedblright
would imply that only countable many measures can be computed whereas, in
general, a dynamical system may have uncountably many invariant measures
(usually an infinite dimensional set). So, a priori most of them will not be
algorithmically describable. This is not a problem because we should put our
attention on the most \textquotedblleft meaningful\textquotedblright\ ones.
An important part of the theory of dynamical systems is indeed devoted to
the understanding of \textquotedblleft physically\textquotedblright\
relevant invariant measures. Informally speaking, these are measures which
represent the asymptotic statistical behavior of \textquotedblleft
many\textquotedblright\ (positive Lebesgue measure) initial conditions (see
Section \ref{DS} for more details).

The existence and uniqueness of physical measures is a widely studied
problem (see \cite{Y}), which has been solved for some important classes of
dynamical systems. These measures are some of the good candidates to be
computed.

The more or less \emph{rigorous }computation of such measures is the main
goal of a part of the literature \ related to computation and dynamics. A
main role here is played by the transfer operator induced by the dynamics.
Indeed, the map $T$ \ defining the dynamics, induces a dynamics $L_{T}$ on
the space of probability measures on $X,$ $L_{T}:PM(X)\rightarrow PM(X)$. $%
L_{T}$ is called the transfer operator associated to $T$ (definition and
basic results about this are recalled in Section \ref{DS}). Invariant
measures are fixed points of this operator. The main part of the methods
which are used to compute invariant measures deals with suitable finite
dimensional approximation of this operator. In Section \ref{DS} we will
review briefly some of these methods, and some references.We then consider
the problem from an abstract point of view, and give some general result on
rigorous computability of the physical invariant measure. In particular we
will see that the transfer operator is computable up to any approximation in
a general context (see Thm \ref{L_comp}). The invariant measure is
computable, provided we are able to give a description of \ a space of
\textquotedblleft regular\textquotedblright\ measures where the physical
invariant measure is the only invariant one (see Thm \ref{comp_inv} and
following corollaries).

We will also show how the description of the space can be obtained from the
one of the system in a class of examples (piecewise expanding maps) by using
the Lasota-Yorke inequality.

After such general statements one could conjecture that all computable
dynamical systems should always have a computable invariant measure. We will
see that, perhaps surprisingly, this is not true. Not all systems that can
be described explicitly (the dynamics can be computed up to any prescribed
approximation) have computable invariant measures (see Section \ref{lst}).
The existence of such examples reveals some subtleties in the computation of
invariant measures.

To further motivate these results, we finally remark that from a technical
point of view, computability of the considered invariant measure is a
requirement in several results about relations between computation,
probability, randomness and pseudo-randomness (see Section \ref{pseudorandom}
and e.g. \cite{LM08, GHR08, GHR07, GHR09c})

\medskip

\noindent\textbf{Computability in Complex Dynamics}

Polynomial Julia sets have emerged as the most studied examples of fractal
sets generated by a dynamical system. One of the reasons for their
popularity is the beauty of the computer-generated images of such sets. The
algorithms used to draw these pictures vary; the na\"ive approach in this
case works by iterating the center of a pixel to determine if it lies in the
Julia set. There exists also more sophisticated algorithms (using classical
complex analysis, see \cite{Mil}) which work quite well for many examples,
but it is well known that in some particular cases computation time will
grow very rapidly with increase of the resolution. Moreover, there are
examples, even in the family of quadratic polynomials, where no satisfactory
pictures of the Julia set exist.

In the rigorous approach, a set is computable if, roughly speaking, its
image can be generated by a computer with an arbitrary precision. Under this
notion of computability, the question arise if Julia sets are always
computable. In a series of papers (\cite{BBY06, BBY07, BY06, BY08}) it was
shown that even though in many cases (hyperbolic, parabolic) the Julia set
is indeed computable, there exists quadratic polynomials which are
computable (again, in the sense that all the trajectories can be
approximated by an algorithm at any desired accuracy), and yet the Julia set
is not.

So we can not simulate the set of limits points on which the chaotic
dynamics takes place, but, what about the statistical distribution?. In
fact, it was shown by Brolin and Luybich that there exists a unique
invariant measure which maximizes entropy, and that this measure is
supported on the Julia set. The question of whether this measure can be
computed has been recently solved in \cite{BBYR11}, where it is proved that
the Brolin-Lyubich measure is always computable. So that even if we can not
visualize the Julia set as a spatial object, we can approximate its
distribution at any finite precision.

\medskip

\noindent\textbf{Computing the speed of convergence and pseudorandom points.}

In several questions in ergodic theory \ The knowledge of the speed of
convergence to ergodic behavior is important to deduce other practical
consequences. In the computational framework, the question turn out to be
the effective estimation of the speed of convergence in the ergodic theorems%
\footnote{%
Find a $N$ such that $\frac{1}{n}\sum f\circ T^{n}$ $\ $differs from $\int \!%
{f}\,\mathrm{d}{\mu }$ less than a given error for each $n\geq N$.}. From
the numerical-practical point of view this has been done in some classes of
systems, having a spectral gap for example. In this case a suitable
approximation of the transfer operator allows to compute the rate of decay
of correlations \cite{FCMP}\cite{L} and from this, other rates of
convergence can be easily deduced.

Other classes of systems could be treated joining the above spectral
approach, with combinatorial constructions (towers, see \cite{LH} e.g.), but
the general case need a different approach.

In \cite{AvigadGT10} it was shown that much more in general, \emph{if the
system can be described effectively, then the rate of convergence in the
pointwise ergodic theorem can be effectively estimated}. We give in section %
\ref{speed} a very short proof of a statement of this kind (see Theorem \ref%
{corollaryx} ) for ergodic systems, and show some consequences. Among these,
a \emph{constructive} version of pointwise ergodic theorem. If the system is
computable (in some wide sense that will be described) then, it is possible
to compute points having typical statistical behavior. Such points could be
hence called \emph{pseudorandom }points in the system (see Section \ref%
{pseudorandom}).

Since the computer can only handle computable initial conditions, any
simulation can start only from these points. Pseudorandom initial conditions
are hence in principle good points where to start a simulation.

We remark that it is widely believed that naive computer simulations very
often produce correct statistical behavior. The evidence is mostly
heuristic. Most arguments are based on the various \textquotedblleft
shadowing\textquotedblright\ results (see e.g. \cite{KH} chapter 18). In
this kind of approach (different from ours), it is possible to prove that in
a suitable system every pseudo-trajectory, as the ones which are obtained in
simulations with some computation error, is close to a real trajectory of
the system. However, even if we know that what we see in a simulation is
near to some real trajectory, we do not know if this real trajectory is
typical in some sense. A limit of this approach is that shadowing results
hold only in particular systems, having some strong hyperbolicity, while
many physically interesting systems are not like this. In our approach we
consider real trajectories instead of "pseudo" ones and we ask if there is
some computable point which behaves as a typical point for the dynamics.

\begin{acknowledgement}
We would like to thank The Abdus Salam International Centre for Theoretical
Physics (Trieste, IT) for support and hospitality during this research.
\end{acknowledgement}

\section{Computability on metric spaces}

To have formal results and precise assumptions on the computability (up to
any given error) of continuous objects, we have to introduce some concepts.

We have to introduce some recursive version of open and compact sets, and
characterize the functions which are well suited to operate with those sets
\ (computable functions). We explain this theory in this section trying the
explanation to be to be as much as possible simple and self contained.

\subsection{Computability}

The starting point of recursion theory was to give a mathematical definition
making precise the intuitive notions of \emph{algorithmic} or \emph{%
effective procedure} on symbolic objects. Several different formalizations
have been independently proposed (by Church, Kleene, Turing, Post,
Markov...) in the 30's, and have proved to be equivalent: they compute the
same functions from $\mathbb{N}$ to $\mathbb{N}$. This class of functions is
now called the class of \emph{recursive functions}. As an algorithm is
allowed to run forever on an input, these functions may be \emph{partial},
i.e. not defined everywhere. The \emph{domain} of a recursive function is
the set of inputs on which the algorithm eventually halts. A recursive
function whose domain is $\mathbb{N}$ is said to be \emph{total}.

We now recall an important concept from recursion theory. A set $E\subseteq 
\mathbb{N}$ is said to be \textbf{\emph{recursively enumerable (r.e.)}} if
there is a (partial or total) recursive function $\varphi :\mathbb{N}%
\rightarrow \mathbb{N}$ enumerating $E$, that is $E=\{\varphi (n):n\in 
\mathbb{N}\}$. If $E\neq \emptyset $, $\varphi $ can be effectively
converted into a total recursive function $\psi $ which enumerates the same
set $E$.

\subsection{Algorithms and uniform algorithms}

Strictly speaking, recursive functions only work on natural numbers, but
this can be extended to the objects (thought as ``finite'' objects) of any
countable set, once a numbering of its elements has been chosen. We will use
the word \emph{algorithm} instead of \emph{recursive function} when the
inputs or outputs are interpreted as finite objects. The operative power of
algorithms on the objects of such a numbered set obviously depends on what
can be effectively recovered from their numbers.

More precisely, let $X$ and $Y$ be numbered sets such that the numbering of $%
X$ is injective (it is then a bijection between $\mathbb{N}$ and $X$). Then
any \emph{recursive function} $\varphi:\mathbb{N}\to\mathbb{N} $ induces an 
\emph{algorithm} $\mathcal{A}:X\to Y$. The particular case $X=\mathbb{N}$
will be much used.

For instance, the set $\mathbb{Q}$ of rational numbers can be injectively
numbered $\mathbb{Q}=\{q_0,q_1,\ldots\}$ in an \emph{effective} way: the
number $i$ of a rational $a/b$ can be computed from $a$ and $b$, and vice
versa. We fix such a numbering: from now and beyond the rational number with
number $i$ will be denoted by $q_i$.

Now, let us consider computability notions on the set $\mathbb{R}$ of real
numbers, introduced by Turing in \cite{T36}.

\begin{definition}
Let $x$ be a real number. We say that:

\begin{itemize}
\item $x$ is \textbf{\emph{lower semi-computable}} if the set $\{i\in 
\mathbb{N}:q_i<x\}$ is r.e.

\item $x$ is \textbf{\emph{upper semi-computable}} if the set $\{i\in 
\mathbb{N}:q_i>x\}$ is r.e.

\item $x$ is \textbf{\emph{computable}} if it is lower and upper
semi-computable.
\end{itemize}
\end{definition}

Equivalently, a real number is computable if and only if there exists an
algorithmic enumeration of a sequence of rational numbers converging
exponentially fast to $x$. That is:

\begin{proposition}
A real number is \textbf{\emph{computable}} if there is an algorithm $%
\mathcal{A}:\mathbb{N}\to \mathbb{Q}$ such that $|\mathcal{A}(n)-x|\leq
2^{-n}$ for all $n$.
\end{proposition}

\medskip

\noindent \textbf{Uniformity.} Algorithms can be used to define
computability notions on many classes of mathematical objects. The precise
definitions will be particular to each class of objects, but they will
always follow the following scheme:

%\begin{center}
%An object $O$ is \textbf{\emph{computable}} if there is an \\[0pt]
%algorithm $\mathcal{A}:X\to Y$ which computes $O$ in some way.
%\end{center}

\smallskip

\begin{center}
An object $O$ is \textbf{\emph{computable}} if there is an algorithm 
\begin{equation*}
\mathcal{A}:X\to Y
\end{equation*}

which computes $O$ in some way.
\end{center}

\smallskip

Each computability notion comes with a uniform version. Let $(O_i)_{i\in%
\mathbb{N}}$ be a sequence of computable objects:

\smallskip

\begin{center}
$O_i$ is computable \textbf{\emph{uniformly in $\boldsymbol{i}$}} if there
is an algorithm 
\begin{equation*}
\mathcal{A}:\mathbb{N}\times X \to Y
\end{equation*}

such that for all $i$, $\mathcal{A}_i:= \mathcal{A}(i,\cdot):X\to Y$
computes $O_i$.
\end{center}

\smallskip

For instance, the elements of a sequence of real numbers $(x_{i})_{i\in 
\mathbb{N}}$ are uniformly computable if there is a algorithm $\mathcal{A}:%
\mathbb{N}\times \mathbb{N}\rightarrow \mathbb{Q}$ such that $|\mathcal{A}%
(i,n)-x_{i}|\leq 2^{-n}$ for all $i,n$.

In other words a set of objects is computable uniformly with respect to some
index if they can be computed with a "general" algorithm \emph{starting}
from the value of the index.

In each particular case, the computability notion may take a particular
name: computable, recursive, effective, r.e., etc. so the term
``computable'' used above shall be replaced. %*****************************

\subsection{Computable metric spaces\label{CMS}}

A computable metric space is a metric space with an additional structure
allowing to interpret input and output of algorithms as points of the metric
space. This is done in the following way: there is a dense subset (called
ideal points) such that each point of the set is identified with a natural
number. The choice of this set is compatible with the metric, in the sense
that the distance between two such points is computable up to any precision
by an algorithm getting the names of the points as input. Using these simple
assumptions many constructions on metric spaces can be implemented by
algorithms.

\begin{definition}
A \textbf{\emph{computable metric space}} (CMS) is a triple $\mathcal{X}%
=(X,d,S)$, where

\begin{enumerate}
\item[(i)] $(X,d)$ is a separable metric space.

\item[(ii)] $S=\{s_{i}\}_{i\in \mathbb{N}}$ is a dense, numbered, subset of $%
X$ called the set of \textbf{\emph{ideal points}}.

\item[(iii)] The distances between ideal points $d(s_{i},s_{j})$ are all
computable, uniformly in $i,j$ (there is an algorithm $\mathcal{A}:\mathbb{N}%
^3\to \mathbb{Q}$ such that $|\mathcal{A}(i,j,n)-d(s_i,s_j)|<2^{-n}$).
\end{enumerate}
\end{definition}

$S$ is a numbered set, and the information that can be recovered from the
numbers of ideal points is their mutual distances. Without loss of
generality, we will suppose the numbering of $S$ to be injective: it can
always be made injective in an effective way.

\begin{definition}
\label{recconv}We say that in a metric space $(X,d)$, a sequence of points $%
(x_{n})_{n\in \mathbb{N}}$ converges\emph{\ recursively }to a point $x$ if
there is an algorithm $D:\mathbb{Q\rightarrow N}$ such that $d(x_{n},x)\leq
\epsilon $ for all $n\geq D(\epsilon )$.
\end{definition}

\begin{definition}
\label{comp_points}A point $x\in X$ is said to be \textbf{\emph{computable}}
if there is an algorithm $\mathcal{A}:\mathbb{N}\rightarrow S$ such that $(%
\mathcal{A}(n))_{n\in \mathbb{N}}$ converges recursively to $x$.
\end{definition}

We define the set of \textbf{\emph{ideal balls}} to be $\mathcal{B}%
:=\{B(s_{i},q_{j}):s_{i}\in S,0<q_{j}\in \mathbb{Q}\}$ where $B(x,r)=\{y\in
X:d(x,y)<r\}$ is an open ball. We fix a numbering $\mathcal{B}%
=\{B_{0},B_{1},\ldots \}$ which makes the number of a ball effectively
computable from its center and radius and vice versa. $\mathcal{B}$ is a
countable basis of the topology.

\begin{definition}[Effective open sets]
\label{reop} We say that an open set $U$ is \emph{effective} if there is an
algorithm $\mathcal{A}:\mathbb{N}\rightarrow \mathcal{B}$ such that $%
U=\bigcup_{n}\mathcal{A}(n)$.
\end{definition}

Observe that an algorithm which diverges on each input $n$ enumerates the
empty set, which is then an effective open set. Sequences of uniformly
effective open sets are naturally defined. Moreover, if $(U_{i})_{i\in 
\mathbb{N}}$ is a sequence of uniformly effective open sets, then $%
\bigcup_{i}U_{i}$ is an effective open set.

\begin{definition}[Effective $G_{\protect\delta }$-set]
An \textbf{\emph{effective $\boldsymbol{G_{\delta }}$-set}} is an
intersection of a sequence of uniformly effective open sets.
\end{definition}

Obviously, an uniform\ intersection of effective $G_{\delta }$-sets is also
an effective $G_{\delta }$-set.

Let $(X,S_{X}=\{s_{1}^{X},s_{2}^{X},...\},d_{X})$ and $(Y,S_{Y}=%
\{s_{1}^{Y},s_{2}^{Y},...\},d_{Y})$ be computable metric spaces. Let also $%
B_{i}^{X}$ and $B_{i}^{Y}$ be enumerations of the ideal balls in $X$ and $Y$%
. A computable function $X\rightarrow Y$ is a function whose behavior can be
computed by an algorithm up to any precision. For this it is sufficient that
the pre-image of each ideal ball can be effectively enumerated by an
algorithm.

\begin{definition}[Computable Functions]
\label{comp_func} A function $T:X\rightarrow Y$ is \textbf{\emph{computable}}
if $T^{-1}(B_{i}^{Y})$ is an effective open set, uniformly in $i$. That is,
there is an algorithm $\mathcal{A}:\mathbb{N}\times \mathbb{N}\rightarrow 
\mathcal{B}^{X}$ such that $T^{-1}(B_{i}^{Y})=\bigcup_{n}\mathcal{A}(i,n)$
for all $i$.

\noindent A function $T:X\rightarrow Y$ is \textbf{\emph{computable on $%
D\subseteq X$}} if there are uniformly effective open sets $U_{i}$ such that 
$T^{-1}(B_{i}^{Y})\cap D=U_{i}\cap D.$
\end{definition}

\begin{remark}
Intuitively, a function $T$ is computable (on some domain $C$) if there is a
computer program which computes $T(x)$ (for $x\in C$) in the following
sense: on input $\epsilon>0$, the program, along its run, asks the user for
approximations of $x$, and eventually halts and outputs an ideal point $s\in
Y$ satisfying $d(T(x),s)<\epsilon$. This idea can be formalized, using for
example the notion of \emph{oracle computation}. The resulting notion
coincides with the one given in the previous definitions.
\end{remark}

Recursive compactness is an assumption which will be needed in the
following. Roughly, a compact set is recursively compact if the fact that it
is covered by a finite collection of ideal balls can be tested
algorithmically (for equivalence with the $\epsilon $-net approach and other
properties of recursively compact set see \cite{GalHoyRoj3}).

\begin{definition}
\label{rec.compact} A set $K\subseteq X$ is \emph{\textbf{recursively compact%
}} if it is compact and there is a recursive function $\varphi :\mathbb{%
N^{\ast }}\rightarrow \mathbb{N}$ such that $\varphi (i_{1},\ldots ,i_{p})$
halts if and only if $(B_{i_{1}},\ldots ,B_{i_{p}})$ is a covering of $K$.
\end{definition}

\section{Computing invariant measures\label{DS}}

\subsection{Invariant measure and statistical properties}

Let $X$ be a metric space, $T:X\mapsto X$ a Borel measurable map and $\mu $
a $T$-invariant Borel probability measure. A set $A$ is called $T$-invariant
if $T^{-1}(A)=A\ (mod\ 0)$. The system $(X,T,\mu )$ is said to be ergodic if
each $T$-invariant set has total or null measure. In such systems the famous
Birkhoff ergodic theorem says that time averages computed along $\mu $%
-typical orbits coincides with space average with respect to $\mu .$ More
precisely, for any $f\in L^{1}(X,\mu )$ it holds 
\begin{equation}
\underset{n\rightarrow \infty }{\lim }\frac{S_{n}^{f}(x)}{n}=\int \!{f}\,%
\mathrm{d}{\mu },  \label{Birkhoff}
\end{equation}%
for $\mu $ almost each $x$, where $S_{n}^{f}=f+f\circ T+\ldots +f\circ
T^{n-1}.$

This shows that in an ergodic system, the statistical behavior of
observables, under typical realizations of the system is given by the
average of the observable made with the invariant measure.

We say that a point $x$ belongs to the basin of an invariant measure $\mu $
if (\ref{Birkhoff}) holds at $x$ for each bounded continuous $f$. In case $X$
is a manifold (possibly with boundary), a physical measure is an invariant
measure whose basin has positive Lebesgue measure (for more details and a
general survey see \cite{Y}). Computation of such measures will be the main
subject of this section.

\subsubsection{The transfer operator\label{PF}}

A function $T$ between metric spaces naturally induces a function $L_{T}$
between probability measure spaces. This function $L_{T}$ is linear and is
called transfer operator (associated to $T$). Measures which are invariant
for $T$ are fixed points of $L_{T}$.

Let us consider a computable metric space $X$ and a measurable function $%
T:X\rightarrow X$. Let us also consider the space $PM(X)$ of Borel
probability measures on $X.$

Let us define the linear function $L_{T}:PM(X)\rightarrow PM(X)$ by duality
in the following way: if $\mu \in PM(X)$ then $L_{T}(\mu )$ is such that

\begin{equation*}
\int \!{f}\,\mathrm{d}{L_{T}(\mu )}=\int \!{f\circ T}\,\mathrm{d}{\mu .}
\end{equation*}

The computation of invariant measures (and many other dynamical quantities)
very often is done by computing the fixed points (or the spectrum) of this
operator in a suitable function space. The most applied and studied strategy
is to find a suitable finite dimensional approximation of $L_{T}$
(restricted to a suitable function space) so reducing the problem to the
computation of the corresponding relevant eigenvectors of a finite matrix.

An example of this is done by discretizing the space $X$ by a partition $%
A_{i}$ and replacing the system by a (finite state) Markov Chain with
transition probabilities%
\begin{equation*}
P_{ij}=\frac{m(A_{i}\cap A_{j})}{m(A_{i})}
\end{equation*}%
where $m$ is the Lebesgue measure on $X$ \ (see e.g. \cite{F07}\cite{F08}%
\cite{L} ), then, taking finer and finer partitions it is possible to obtain
in some cases that the finite dimensional model will converge to the real
one (and its natural invariant measure to the physical measure of the
original system). In some case there is an estimation for this speed of
convergence (see eg. \cite{F08} for a discussion), but a rigorous bound on
the error (and then a real rigorous computation) is known only in a few
cases (piecewise expanding or expanding maps, see \cite{L}).

Similar approaches consists in applying a kind of Faedo-Galerkin \
approximation to the transfer operator by considering a complete Hilbert
base of the function space and truncating the operator to the action on the
first elements (see \cite{web} ).

Another approach is to consider a perturbation of the system by a small
noise. The resulting transfer operator has a kernel. This operator then can
be approximated by a finite dimensional one, again by Faedo-Galerkin \
method and relevant eigenvectors are calculated (see e.g. \cite{Del99}\cite%
{DelJu02}) then, if we prove that the physical measure of the original
system can be obtained as a limit when the size of the noise tends to zero
(this happen on uniformly hyperbolic system for example) we have a method
which in principle can rigorously compute this measure.

Variations on the method of partitions are given in \cite{Din93, Din94},
while in \cite{PJ99} a different method, fastly converging, based on
periodic points is exploited for piecewise analytic Markov maps. Another
strategy to face the problem of computation of invariant measures consist in
following the way the measure $\mu $ can be constructed and check that each
step can be realized in an effective way. In some interesting examples we
can obtain the physical measure as limit of iterates of the Lebesgue measure 
$\mu =\lim_{n\rightarrow \infty }L_{T}^{n}(m)$ (recall that $m$ is the
Lebesgue measure). To prove computability of $\mu $ the main point is to
explicitly estimate the speed of convergence to the limit. This sometimes
can be done using the spectral properties of the system (\cite{GHR07}).

Concluding, if the goal is to rigorously compute an invariant measure, most
of the results which are in the today literature are partial. Indeed, beside
being applied to a quite restricted class of systems, to compute the measure
those methods need additional information. For example, the way to compute
rigorously the finite dimensional approximation is often not done
effectively, or the rate of convergence of the approximation is computed up
to some constants depending on the system, which are not estimated.

In the remaining part of the section we present some results, mainly from 
\cite{GalHoyRoj3} explaining some result about rigorous computations of
invariant measures. These results have the advantage to give in principle an
effective method for the rigorous computation of an invariant measure and
give a quite general result, where all the needed assumptions are explicited.

They have the disadvantage to not optimize computation time. So it is not
clear if they can be implemented and used in practice.

The rigorous framework into they are proved, however allows to see them as
an investigation about the theoretical limits of \ rigorous computation of
invariant measures.

Moreover by this, also negative results can be proved. In particular, we
give examples of \emph{computable }systems having no computable invariant
measures.

\subsection{Computability of measures\label{meas}}

In this section we explain precisely what we mean by computing a measure.
This means, having an algorithm who is able to approximate the measure by
"simple measures" up to any given approximation.

Let us consider the space $PM(X)$ of Borel probability measures over $X$.
Let $C_{0}(X)$ be the set of real-valued bounded continuous functions on $X$%
. We recall the notion of weak convergence of measures:

\begin{definition}
$\mu _{n}$ is said to be \emph{\textbf{weakly convergent}} to $\mu $ if $%
\int\!{f}\,\mathrm{d}{\mu_n}\rightarrow \int\!{f}\,\mathrm{d}{\mu}$ for each 
$f\in C_{0}(X)$.
\end{definition}

Let us introduce the Wasserstein-Kantorovich distance between measures. Let $%
\mu _{1}$ and $\mu _{2}$ be two probability measures on $X$ and consider:

\begin{equation*}
W_{1}(\mu _{1},\mu _{2})=\underset{f\in 1\text{-Lip}(X)}{\sup }\left|\int\!{f%
}\,\mathrm{d}{\mu_1}-\int\!{f}\,\mathrm{d}{\mu _2}\right|
\end{equation*}

where $1\mbox{-Lip}(X)$ is the space of functions on $X$ having Lipschitz
constant less than one. The distance $W_{1}$ has the following useful
properties

\begin{proposition}[see \protect\cite{AGS} Prop 7.1.5]
\label{ambros}\mbox{}

\begin{enumerate}
\item $W_{1}$ is a distance and if $X$ is bounded, separable and complete,
then $PM(X) $ with this distance is a separable and complete metric space.

\item If $X$ is bounded, a sequence is convergent for the $W_{1}$ metrics if
and only if it is convergent for the weak topology.

\item If $X$ is compact $PM(X)$ is compact with this metrics.
\end{enumerate}
\end{proposition}

Item (1) has an effective version: $PM(X)$ inherits the computable metric
structure of $X$. Indeed, given the set $\S _{X}$ of ideal points of $X$ we
can naturally define a set of ideal points $\S _{PM(X)}$ in $PM(X)$ by
considering finite rational convex combinations of the Dirac measures $%
\delta _{s}$ supported on ideal points $s\in S_{X}$. This is a dense subset
of $PM(X)$. The proof of the following proposition can be found in (\cite%
{HR07}).

\begin{proposition}
If $X$ bounded then $(PM(X),W_1, \S _{PM(X)})$ is a computable metric space.
\end{proposition}

A measure $\mu $ is then computable if there is a sequence $\mu _{n}\in \S %
_{PM(X)}$ converging recursively fast to $\mu $ in the $W_{1}$ metric (and
hence for the weak convergence).

\subsection{Computable invariant \textquotedblleft
regular\textquotedblright\ measures}

\label{section.measures}

Invariant measures can be found as \ fixed points of the transfer operator,
i.e. as solutions of the equation $W_{1}(\mu ,L_{T}(\mu ))=0.$ There is an
abstract theorem allowing the computation of isolated fixed points (see \cite%
{GalHoyRoj3} Corollary 3) of computable maps. To apply it is necessary to
consider a space where the transfer transfer operator is computable (this
space hopefully will contain the physical invariant measure).

We remark that if $T$ is not continuous then $L_{T}$ is not necessarily
continuous (this can be realized by applying $L_{T}$ to some delta measure
placed near a discontinuity point) hence not computable. Still, we have that 
$L_{T}$ is continuous (and its modulus of continuity is computable) at all
measures $\mu $ which are \textquotedblleft far enough\textquotedblright\
from the discontinuity set $D$. This is technically expressed by the
condition $\mu (D)=0$.

\begin{proposition}
\label{L_comp} Let $X$ be a computable metric space and $T:X\rightarrow X$
be a function which is computable on $X\setminus D$. Then $L_{T}$ is
computable on the set of measures 
\begin{equation}
PM_{D}(X):=\{\mu \in PM(X):\mu (D)=0\}.
\end{equation}
\end{proposition}

From a practical point of view this proposition provides sufficient
conditions to rigorously approximate the transfer operator by an algorithm.

The above tools allow us to ensure the computability of $L_{T}$ on a large
class of measures. This will enable us to apply a general result on
computation of fixed points and obtain

\begin{theorem}[\protect\cite{GalHoyRoj3}, Theorem 3.2]
\label{comp_inv} Let $X$ be a computable metric space and $T$ a function
which is computable on $X\setminus D$. Suppose there is a recursively
compact set of probability measures $V\subset PM(X)$ such that for every $%
\mu \in V$, $\mu (D)=0$ holds. Then every invariant measure isolated in $V$
is computable. Moreover the theorem is uniform: there is an algorithm which
takes as inputs finite descriptions of $T,V$ and an ideal ball in $M(X)$
which isolates\footnote{%
If the invariant measure is unique in $V$ the isolating ball is not
necessary.} an invariant measure $\mu $, and outputs a finite description of 
$\mu $.
\end{theorem}

A trivial and general consequence of Theorem \ref{comp_inv} is the following:

\begin{corollary}
If a system as above is uniquely ergodic and its invariant measure $\mu$
satisfies $\mu (D)=0$, then it is a computable measure.
\end{corollary}

The main problem in the application of Theorem \ref{comp_inv} is the
requirement that the invariant measure we are trying to compute be isolated
in $V$. In general the space of invariant measures in a given dynamical
system could be very large (an infinite dimensional convex subset of $PM(X)$%
) and physical measures have often some kind of particular regularity. To
isolate a particular measure we can restrict and consider a subclass of
\textquotedblleft regular\textquotedblright\ measures.

Let us consider the following \emph{seminorm}:

\begin{equation*}
\left\Vert \mu \right\Vert _{\alpha }=\sup_{x\in X,r>0}\frac{\mu (B(x,r))}{%
r^{\alpha }}.
\end{equation*}%
If $\alpha $ and $K$ are computable and $X$ is recursively compact then 
\begin{equation}
V_{\alpha ,K}=\{\mu \in PM(X):\left\Vert \mu \right\Vert _{\alpha }\leq K\}
\end{equation}%
is recursively compact (\cite{GalHoyRoj3}). This implies

\begin{proposition}
\label{2}Let $X$ be recursively compact and $T$ be computable on $X\setminus
D$, with $dim_{H}(D)<\infty $. Then any invariant measure isolated in $%
V_{\alpha ,K}$ with $\alpha >dim_{H}(D)$ is computable. Once again, this is
uniform in $T,\alpha ,K$.
\end{proposition}

The above general propositions allow us to obtain as a corollary the
computability of many absolutely continuous invariant measures. For the sake
of simplicity, let us consider maps on the interval.

\begin{proposition}
\label{prop_dim} If $X=[0,1]$, $T$ is computable on $X\setminus D$, with $%
dim_{H}(D)<1$ and $(X,T)$ has a unique absolutely continuous invariant
measure $\mu $ with bounded density, then $\mu $ is computable (starting
from $T$ and a bound for the $L^{\infty }$ norm of the invariant density).
\end{proposition}

Similar results hold for maps on manifolds (again see \cite{GalHoyRoj3}).

\subsubsection{Computing the measure from a description of the system in the
class of piecewise expanding maps\label{PWSEC}}

As it is well known, interesting examples of systems having a unique
absolutely continuous invariant measure (with bounded density as required)
are topologically transitive \emph{piecewise expanding maps} on the interval
or \emph{expanding maps} on manifolds.

We show how to find a bound for the invariant density on piecewise expanding
maps. By this, the invariant measure can be calculated starting from a
description of $T$.

\begin{definition}
A nonsingular function $T:([0,1],m)\rightarrow ([0,1],m)$ is said piecewise
expanding if\footnote{%
For the seek of simplicity we will consider the simplest setting in which we
can work and give precise estimations. Such class was generalized in several
ways, we hence warn the rader that now in the current literature by
piecewise expanding maps it it meant something slightly more general than
the ones defined here.}

\begin{enumerate}
\item There is a finite set of points $d_{1}=0,d_{2},...,d_{n}=1$ such that $%
T|_{(d_{i},d_{i+1})}$ is $C^{2}$ and can be extended to a $C^{2}$map on $%
[d_{i},d_{i+1}]$.

\item $\inf (T^{\prime })>1$ on the set where it is defined.

\item $T$ is topologically mixing.
\end{enumerate}
\end{definition}

It is now well known that such maps have an unique ergodic invariant measure
with bounded variation density.

Such density is also the unique fixed point of the (Perron Frobenius )
operator\footnote{%
Note that this operator corresponds to the above cited transfes operator,
but it acts on densities instead of measures.} $L:L^{1}[0,1]\rightarrow
L^{1}[0,1]$ defined by

\begin{equation*}
\lbrack Lf](x)=\sum_{y\in T^{-1}x}\frac{f(y)}{T^{\prime }(y)}.
\end{equation*}

We now show how to find a bound for such density, starting from the
description of the system, and then compute the associated invariant measure.

The following proposition was proved in the celebrated paper\ \cite{LY}
(Thm. 1 and its proof) and it is now called Lasota-Yorke inequality. We give
a precise statement where the constants are explicited.

\begin{proposition}
\label{proe}Let $T$ a piecewise hyperbolic map. If $f$ is of bounded
variation in $[0,1]$. Let $d_{1},...,d_{n}$ be the discontinuity points of $%
T $.

If $\lambda =\underset{x\in \lbrack 0,1]-\{d_{1},...,d_{n}\}}{\inf }%
T^{\prime }(x)$. Then%
\begin{equation*}
Var(Lf)\leq 2\lambda Var(f)+B||f||_{1}
\end{equation*}%
where 
\begin{equation*}
B=\frac{\underset{x\in \lbrack 0,1]-\{d_{1},...,d_{n}\}}{\sup }(|\frac{%
T^{\prime \prime }(x)}{(T^{\prime }(x))^{2}}|)}{\underset{x\in \lbrack
0,1]-\{d_{1},...,d_{n}\}}{\inf }(\frac{1}{T^{\prime }(x)})}+\frac{2}{\min
(d_{i}-d_{i+1})}.
\end{equation*}
\end{proposition}

The following is an elementar fact about the behavior of real sequences

\begin{lemma}
\label{sequ}If a real sequence $a_{n}$ is such that $a_{n+1}\leq la_{n}+k$ \
for some $l<1,k>0$, then%
\begin{equation*}
\sup (a_{n})\leq \max (a_{0},\frac{k}{1-l}).
\end{equation*}
\end{lemma}

\begin{proposition}
\label{densibound}If $f$ is the density of the physical measure of \ a
piecewise expanding map $T$ as above and $\lambda >2$. Then%
\begin{equation*}
Var(f)\leq \frac{B}{1-2\lambda }
\end{equation*}%
where $B$ is defined as above.
\end{proposition}

\begin{proof}
(sketch) The topological mixing assumption implies that the map has only one
invariant physical measures (see \cite{V97}). Let us use the above results
iterating the constant density corresponding to the Lebesgue measure.
Proposition \ref{proe} and Lemma \ref{sequ} give that the variation of the
iterates is bounded by $\frac{B}{1-2\lambda }$. Suppose that the limit
measure has density $f$. By compactness of $BV$ in $L^{1}$, the above
properties give a bound on the variation of $f$ (see \cite{LY} proof of Thm.
1).
\end{proof}

The following is a trivial consequence of the fact that $||f||_{\infty }\leq
Var(f)+\int fd\mu $.

\begin{corollary}
In the above situation $||f||_{\infty }\leq \frac{B}{1-2\lambda }+1$.
\end{corollary}

The bound on the density of the invariant measure, together with Corollary %
\ref{prop_dim} gives the following \emph{uniform} result on the computation
of invariant measures of such maps.

\begin{theorem}
\label{pwexpmu}Suppose a piecewise expanding map $T$ and also its derivative
is computable on $[0,1]-\{d_{1},...,d_{n}\}$, and also its extensions to the
closed intervals $[d_{i},d_{i+1}]$ are computable. Then, the physical
measure can be computed starting from a description of the system (the
program computing the map and its derivative).
\end{theorem}

\begin{proof}
(sketch) Since $T$ \ and $T^{\prime }$ are computable on each interval $%
[d_{i},d_{i+1}]$ we can compute a number $\lambda $ such that $1<\lambda
\leq \inf (T^{\prime })$ (see \cite{GalHoyRoj3} Proposition 3)\ \ If we
consider the iterate $T^{N}$instead of $T$ the associated invariant density
will be the same as the one of $T$, and if $\lambda ^{N}>2$,then $T^{N}$
will satisfy all the assumptions needed on Proposition \ref{densibound}.
Then we have a bound on the density and we can apply Corollary \ref{prop_dim}
to compute the invariant measure.
\end{proof}

\subsection{Unbounded densities, non uniformly hyperbolic maps}

The above results ensure computability of some absolutely continuous
invariant measure with bounded density. If we are interested in situations
where the density is unbounded, we can consider a new norm, ``killing''
singularities.

Let us hence consider a computable function $f:X\rightarrow \mathbb{R}$ \ and

\begin{equation*}
\left\Vert \mu \right\Vert _{f,\alpha }=\sup_{x\in X,r>0}\frac{f(x)\mu
(B(x,r))}{r^{\alpha }}.
\end{equation*}%
If $\alpha $ and $K$ are computable and $X$ is recursively compact then 
\begin{equation}
V_{\alpha ,K}=\{\mu \in PM(X):\left\Vert \mu \right\Vert _{f,\alpha }\leq K\}
\end{equation}%
is recursively compact, and \ref{2} also hold for this norm. If $f$ is such
that $f(x)=0$ when $\lim_{r\rightarrow 0}\frac{\mu (B(x,r))}{r^{\alpha }}%
=\infty $ this can let the norm be finite when the density diverges.

As an example, where this can be applied, let us consider the Manneville
Pomeau type maps on the unit interval. These are maps of the type $%
x\rightarrow x+x^{z}(\mathop{\rm mod}\nolimits~1)$.

When $1<z<2$ the dynamics has a unique absolutely continuous invariant
measure $\mu _{z}$. This measure has a density $e_{z}(x)$ which diverges at
the origin, and $e_{z}(x)\asymp x^{-z+1}$ and is bounded elsewhere (see \cite%
{I03} Section 10 and \cite{V97} Section 3 e.g.). If we consider the norm $%
\left\Vert .\right\Vert _{f,1}$ with $f(x)=x^{2}$ we have that $\left\Vert
\mu _{z}\right\Vert _{f,1}$ is finite for each such $z$. By this it follows
that for each such $z$ the measure $\mu _{z}$ is computable.

\subsection{Computable systems without computable invariant measures\label%
{lst}}

We have seen several techniques to establish the computability of many
physical invariant measures. This raises naturally the following question: a
computable systems does necessarily have a computable invariant measure?
what about ergodic physical measures?

The following is an easy example showing that this is not true in general
even in quite simple systems, hence the whole question of computing
invariant measures has some subtlety.

Let us consider a system on the unit interval given as follows. Let $\tau
\in (0,1)$ be a lower semi-computable real number which is not computable.
There is a computable sequence of rational numbers $\tau _{i}$ such that $%
\sup_{i}\tau _{i}=\tau $. For each $i$ consider $T_{i}(x)=\max (x,\tau _{i})$
and 
\begin{equation*}
T(x)=\sum_{i\geq 1}2^{-i}T_{i}.
\end{equation*}
The functions $T_{i}$ are uniformly computable so $T$ is also computable.

\begin{figure}[ht]
\begin{center}
\includegraphics[height=5cm]{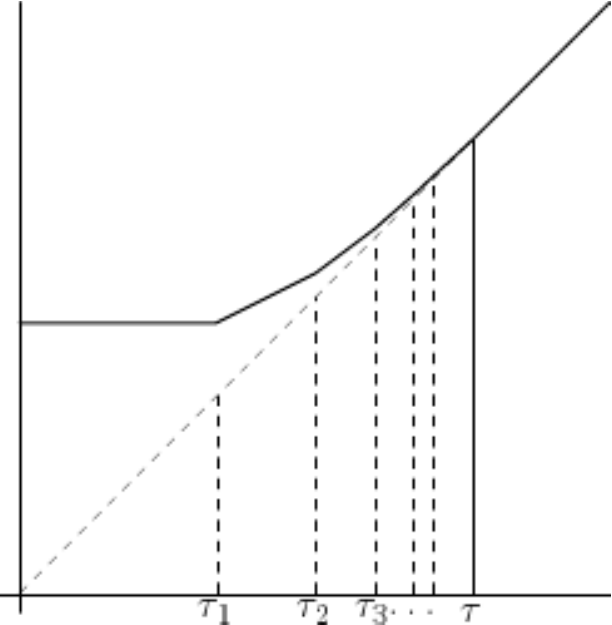}
\end{center}
\caption{The map $T$.}
\end{figure}

The system $([0,1],T)$ is hence a computable dynamical system. $T$ is
non-decreasing, and $T(x)>x$ if and only if $x<\tau $.

This system has a physical ergodic invariant measure which is $\delta _{\tau
}$, the Dirac measure placed on $\tau $. The measure is physical because $%
\tau $ attracts all the interval at its left. Since $\tau $ is not
computable then $\delta _{\tau }$ is not computable. We remark that
coherently with the previous theorems $\delta _{\tau }$ is not isolated.

It is easy to prove, by a simple dichotomy argument, that a computable
function from $[0,1]$ to itself must have a computable fixed point. Hence it
is not possible to construct a system over the interval having no computable
invariant measure (we always have the $\delta $ over the fixed point). With
some more work we will see that such an example can be constructed on the
circle. Indeed the following can be established

\begin{proposition}
\label{bomb} There is a computable, continuous map $T$ on the circle having
no computable invariant probability measure.
\end{proposition}

For the description of the system and for applications to reverse
mathematics see \cite{GalHoyRoj3} (Proposition 12).

\section{Computability in Complex Dynamics}

Let $K$ be a compact subset of the plane. Informally speaking, in order to
draw the set $K$ on a computer screen we need a program which is able to
decide, given some precision $\epsilon $, if pixel $p$ has to be colored or
not. By representing pixel $p$ by a ball $B(p,\epsilon )$ where $(p,\epsilon
)\in \mathbb{Q}^{2}$, the question one would like to answer is: does $%
B(p,\varepsilon )$ intersects $K$ ? The following definition captures
exactly this idea:

\begin{definition}
A compact set $K\subset \mathbb{R}^{2}$ is said to be \textbf{\textit{%
computable}} if there is an algorithm $\mathcal{A}$ such that, upon input $%
(p,\epsilon)$:

\begin{itemize}
\item halts and outputs ``\emph{yes}'' if $K\cap B(p,\epsilon)\neq \emptyset$%
,

\item halts and outputs ``\emph{no}'' if $K\cap \overline{B(p,\epsilon)}%
=\emptyset$,

\item run forever otherwise.
\end{itemize}
\end{definition}

We remark that if a compact set is computable under the definition above,
then it is recursively compact in the sense of Definition \ref{rec.compact}.
The converse is however false.

\subsection{Julia sets and the Brolin-Lyubich measure}

The question of whether Julia sets are computable under this definition has
been completely solved in a series of papers by Binder, Braverman and
Yampolsky. See \cite{BBY07, BBY06, BY06, BY08, BY08a}. Here we review some
of their results. For simplicity, we give the definition of the Julia set
only for quadratic polynomials. Consider a quadratic polynomial 
\begin{equation*}
P_{c}(z)=z^{2}+c : {\mathbb{C}}\to{\mathbb{C}}.
\end{equation*}

Obviously, there exists a number $M$ such that if $|z|>M$, then the iterates
of $z$ under $P$ will uniformly scape to $\infty$. The \textbf{\textit{%
filled Julia set}} is then defined by:

\begin{equation*}
K_{c}=\{z\in{\mathbb{C}}\;\sup_n|P^n(z)|<\infty\}.
\end{equation*}

That is, the set of points whose orbit remains bounded. For the filled Julia
set one has the following result:

\begin{theorem}[\protect\cite{BY08a}]
The filled Julia set is always computable.
\end{theorem}

The \textbf{\textit{Julia set}} can be now defined by:

\begin{equation*}
J_{c}=\partial K_{c},
\end{equation*}
where $\partial(A)$ denotes the boundary of $A$. The Julia set is the 
\textit{repeller} of the dynamical system generated by $P_{c}$. For all but
finitely many points, the limit of the $n$-th preimages $P_{c}^{-n}(z)$
coincides with the Julia set $J_{c}$. The dynamics of $P_{c}$ on the set $%
J_{c}$ is chaotic, again rendering numerical simulation of individual orbits
impractical. Yet Julia sets are among the most drawn mathematical objects,
and countless programs have been written for visualizing them.

In spite of this, the following result was shown in \cite{BY06}.

\begin{theorem}
There exist \emph{computable} quadratic polynomials $P_c(z)=z^2+c$ such that
the Julia set $J_{c}$ is \emph{not computable.}
\end{theorem}

This phenomenon of \textit{non-computability} is rather subtle and rare. For
a detailed exposition, the reader is referred to the monograph \cite{BY08}.

Thus, we cannot accurately simulate the set of limit points of the preimages 
$(P_c)^{-n}(z)$, but what about their statistical distribution? The question
makes sense, as for all $z\neq \infty$ and every continuous test function $%
\psi$, the averages 
\begin{equation*}
\frac{1}{2^n}\sum_{w\in (P_c)^{-n}(z)} \psi(w)\underset{n\to\infty}{%
\longrightarrow}\int \psi d\lambda,
\end{equation*}
where $\lambda$ is the \textit{Brolin-Lyubich probability measure} \cite%
{B65, L82} supported on the Julia set $J_{c}$. We can thus ask whether the
the Brolin-Lyubich measure is computable. Even if $J_{c}=\func{Supp}%
(\lambda) $ is not a computable set, the answer does not \textit{a priori}
have to be negative. In fact, the following result holds:

\begin{theorem}[\protect\cite{BBYR11}]
The Brolin-Lyubich measure is always computable.
\end{theorem}

The proof of the previous result does not involve (perhaps surprisingly)
much analytic machinery, but it follows from some general principles in the
same spirit of the ideas introduced in Section \ref{section.measures}.

\subsection{The Mandelbrot set}

Another important object in complex dynamics is the Mandelbrot set, which is
defined as follows:

\begin{equation*}
\mathcal{M}:=\{c\in {\mathbb{C}} : 0 \in K_{c} \}.
\end{equation*}

That is, the set of parameters of quadratic polynomials for which the orbit
of $0$ remains bounded. As for Julia sets, there exists many computer
programs to visualize it, and the question arise whether this set is
actually computable in a rigorous sense. This is, up to date, an open
question. However, some partial results have been obtained. For instance, in 
\cite{Her05} it is proved that the Mandelbrot set is recursively compact in
the sense of Definition \ref{rec.compact}. There exists also analytical
questions about the Mandelbrot set that remain unsolved. For instance, it is
unknown whether it is \emph{locally connected}. It has been conjectured that
this is indeed the case. Computability of the Mandelbrot set is closely
related to this conjecture, as it has been observed in \cite{Her05}:

\begin{theorem}[\protect\cite{Her05}]
If the Mandelbrot set is locally connected, then it is computable.
\end{theorem}

It is known that local connectivity implies another well known conjecture,
namely the \emph{hyperbolicity conjecture}, which in turn also implies
computability. For a detailed exposition in this subject the reader is
referred to \cite{CG93}.

\section{Computing the speed of ergodic convergence\label{speed}}

As recalled before, the Birkhoff ergodic theorem tells that, if the system
is ergodic, there is a full measure set of points for which the averages of
the values of the observable $f$ along its trajectory (time averages)
coincides with the spatial average of the observable $f$. Similar results
can be obtained for the convergence in the $L^{2}$ norm, and others. \textbf{%
\ }Many, more refined results are linked to the speed of convergence of this
limit. And the question naturally arise, if there is a possibility to
compute this speed of convergence in a sense similar to Definition \ref%
{recconv}.

In the paper \cite{AvigadGT10} some abstract results imply that in a
computable ergodic dynamical system, the speed of convergence of such
averages can be algorithmically estimated. On the other hand it is also
shown that there are non ergodic systems where this kind of estimations are
not possible. In \cite{GHRP} a very short proof of this result for ergodic
systems is shown. We expose the precise result (Theorem \ref{corollaryx})
and the short proof in this section.

\subsubsection{Convergence of random variables}

We first precise what is meant by "compute the speed" in some pointwise a.e.
convergence.

\begin{definition}
A \textbf{\emph{random variable}} on $(X,\mu)$ is a measurable function $%
f:X\to \mathbb{R}$.
\end{definition}

\begin{definition}
Random variables $f_n$ \textbf{\emph{effectively converge in probability}}
to $f$ if for each $\epsilon>0$, $\mu\{x: |f_n(x)-f(x)|<\epsilon\}$
converges effectively to $1$, uniformly in $\epsilon$. That is, there is a
computable function $n(\epsilon,\delta)$ such that for all $n\geq
n(\epsilon,\delta)$, $\mu\{|f_n-f|\geq\epsilon\}<\delta$.
\end{definition}

\begin{definition}
\label{def_rec_conv}Random variables $f_{n}$ \textbf{\emph{effectively
converge almost surely}} to $f$ if $f_{n}^{\prime }=\sup_{k\geq n}|f_{k}-f|$
effectively converge in probability to $0$.
\end{definition}

\begin{definition}
\label{a.e.compufunct}A \textbf{\emph{computable probability space}} is a
pair $(X,\mu )$ where $X$ is a computable metric space and $\mu $ a
computable Borel probability measure on $X$.
\end{definition}

Let $Y$ be a computable metric space. A function $(X,\mu )\rightarrow Y$ is 
\textbf{\emph{almost everywhere computable}} (a.e. computable for short) if
it is computable on an effective $G_{\delta }$-set of measure one, denoted
by $\mathrm{dom}f$ and called the \emph{domain of computability of $f$}.

A \textbf{\emph{morphism}} of computable probability spaces $f:(X,\mu
)\rightarrow (Y,\nu )$ is a morphism of probability spaces which is a.e.
computable.

\begin{remark}
\label{com_sum} A sequence of functions $f_{n}$ is uniformly a.e. computable
if the functions are uniformly computable on their respective domains, which
are uniformly effective $G_{\delta }$-sets. Observe that in this case,
intersecting all the domains provides an effective $G_{\delta }$-set on
which all $f_{n}$ are computable. In the following we will apply this
principle to the iterates $f_{n}=T^{n}$ of an a.e. computable function $%
T:X\rightarrow X$, which are uniformly a.e. computable.
\end{remark}

\begin{remark}
\label{l1comp} The space $L^{1}(X,\mu )$ (resp. $L^{2}(X,\mu )$) can be made
a computable metric space, choosing some dense set of bounded computable
functions as ideal elements. We say that an integrable function $%
f:X\rightarrow \overline{\mathbb{R}}$ is $L^{1}(X,\mu )$-computable if its
equivalence class is a computable element of the computable metric space $%
L^{1}(X,\mu )$. Of course, if $f=g$ $\mu $-a.e., then $f$ is $L^{1}(X,\mu )$%
-computable if and only if $g$ is. Basic operations on $L^{1}(X,\mu )$, such
as addition, multiplication by a scalar, $min$, $max$ etc. are computable.
Moreover, if $T:X\rightarrow X$ preserves $\mu $ and $T$ is a.e. computable,
then $f\rightarrow f\circ T$ (from $L^{1}$ to $L^{1}$) is computable (see 
\cite{HoyRojCiE09}).\bigskip
\end{remark}

Let us call $(X,\mu ,T)$ a \textbf{\emph{computable ergodic system}} if $%
(X,\mu )$ is a computable probability space where $T$ is a measure
preserving morphism and $(X,\mu ,T)$ is ergodic. Let $||f||$ denote the $%
L^{1}$ norm or the $L^{2}$ norm.

\begin{proposition}
Let $(X,\mu ,T)$ be a computable ergodic system. Let $f$ be a computable
element of $L^{1}(X,\mu )$ (resp. $L^{2}(X,\mu )$).

The $L^{1}$ convergence (resp. $L^{2}$ convergence) of the Birkhoff averages 
$A_{n}=(f+f\circ T+\ldots +f\circ T^{n-1})/n$ is effective. That is: there
is an algorithm $n:\mathbb{Q\rightarrow N}$ such that for each $m\geq
n(\epsilon )$ $||A_{m}-\int fd\mu ||\leq \epsilon $. Moreover the algorithm
depends effectively on $T,\mu ,f$.
\end{proposition}

%\begin{proof}
\textbf{Proof.} Replacing $f$ with $f-\int fd\mu $, we can assume that $\int
fd\mu =0$. The sequence $||A_{n}||$ is computable ( see Remark \ref{l1comp}
) and converges to $0$ by the ergodic theorems.

Given $p\in N$, we write $m\in N$ as $m=np+k$ with $0\leq k<p$. Then

\begin{eqnarray*}
A_{np+k} &=&\frac{1}{np+k}\left( \sum_{i=0}^{n-1}pA_{p}\circ
T^{pi}+kA_{k}\circ T^{pn}\right) \\
||A_{np+k}|| &\leq &\frac{1}{np+k}(np||A_{p}||+k||A_{k}||) \\
&\leq &||A_{p}||+\frac{||{A_{k}}||}{n} \\
&\leq &||{A_{p}}||+\frac{||{f}||}{n}.
\end{eqnarray*}

Let $\epsilon >0$. We can compute some $p=p(\epsilon )$ such that $||{A_{p}}%
||<\epsilon /2$. Then we can compute some $n(\epsilon )\geq \frac{2}{%
\epsilon }||{f}||$. The function $m(\epsilon ):=n(\epsilon )p(\epsilon )$ is
computable and for all $m\geq m(\epsilon )$, $||{A_{m}}||\leq \epsilon $. 
%\end{proof}
$\Box$

\subsection{Effective almost sure convergence}

Now we use the above result to find a computable estimation for the a.s.
speed of convergence.

\begin{theorem}
\label{corollaryx}Let $(X,\mu ,T)$ be a computable ergodic system. If $f$ is 
$L^{1}(X,\mu )$-computable, then the \ a.s. convergence is effective.
Moreover, the rate of convergence can be computed as above starting from $%
T,\mu ,f$.
\end{theorem}

This will be proved by the following

\begin{proposition}
\label{theorem_ergodic_effective_as_convergence} If $f$ is $L^{1}(X,\mu )$%
-computable as above, and $||{f}||_{\infty }$ is bounded, then the
almost-sure convergence is effective (uniformly in $f$ and a bound on $||{f}%
||_{\infty }$ and on $T$, $\mu $ ).
\end{proposition}

To prove this we will use the Maximal ergodic theorem:

\begin{lemma}[Maximal ergodic theorem]
\label{lemma_maximal} For $f\in L^{1}(X,\mu )$ and $\delta >0$, 
\begin{equation*}
\mu (\{\sup_{n}|A_{n}^{f}|>\delta \})\leq \frac{1}{\delta }||{f}||_{1}.
\end{equation*}
\end{lemma}

The idea is simple: compute some $p$ such that $||A_p^f||_{1}$ is small,
apply the maximal ergodic theorem to $g:=A_{p}^{f}$, and then there is $%
n_{0} $, that can be computed, such that $A_{n}^{f}$ is close to $A_{n}^{g}$
for $n\geq n_{0}$.

%\begin{proof}
\textbf{Proof }(of Proposition \ref{theorem_ergodic_effective_as_convergence}%
) Let $\epsilon ,\delta >0$. Compute $p$ such that $||{A_{p}^{f}}||\leq
\delta \epsilon /2$. Applying the maximal ergodic theorem to $g:=A_{p}^{f}$
gives: 
\begin{equation}
\mu (\{\sup_{n}|A_{n}^{g}|>\delta /2\})\leq \epsilon .  \label{bound_g}
\end{equation}

Now, $A_{n}^{g}$ is not far from $A_{n}^{f}$: expanding $A_{n}^{g}$, one can
check that 
\begin{equation*}
A_{n}^{g}=A_{n}^{f}+\frac{u\circ T^{n}-u}{np},
\end{equation*}%
where $u=(p-1)f+(p-2)f\circ T+\ldots +f\circ T^{p-2}$. $||{u}||_{\infty
}\leq \frac{p(p-1)}{2}||{f}||_{\infty }$ so if $n\geq n_{0}\geq 4(p-1)||{f}%
||_{\infty }/\delta $, then $||{A_{n}^{g}-A_{n}^{f}}||_{\infty }\leq \delta
/2$. As a result, if $|A_{n}^{f}(x)|>\delta $ for some $n\geq n_{0}$, then $%
|A_{n}^{g}(x)|>\delta /2$. From (\ref{bound_g}), we then derive 
\begin{equation*}
\mu (\{\sup_{n\geq n_{0}}|A_{n}^{f}|>\delta \})\leq \epsilon .
\end{equation*}

As $n_{0}$ can be computed from $\delta $ and $\epsilon $, we get the
result. %\end{proof}
$\Box $

\begin{remark}
\label{remarkboundunif}This result applies uniformly to a uniform sequence
of computable $L^{\infty }(X,\mu )$ observables $f_{n}$.
\end{remark}

We now extend this to $L^{1}(X,\mu )$-computable functions, using the
density of $L^{\infty }(X,\mu )$ in $L^{1}(X,\mu )$.

%\begin{proof}
\textbf{Proof.} (of Theorem \ref{corollaryx}) Let $\epsilon ,\delta >0$. For 
$M\in \mathbb{N} $, let us consider $f_{M}^{\prime }\in L^{\infty }(X,\mu )$
defined as 
\begin{equation*}
f_{M}^{\prime }(x)=\left\{ 
\begin{array}{cc}
\min (f,M) & if~f(x)\geq 0 \\ 
\max (f,-M) & if~f(x)\leq 0.%
\end{array}%
\right.
\end{equation*}

Compute $M$ such that $||{f-f_{M}^{\prime }}||_{1}\leq \delta \epsilon $.
Applying Proposition \ref{theorem_ergodic_effective_as_convergence} to $%
f_{M}^{\prime }$ gives some $n_{0}$ such that $\mu (\{\sup_{n\geq
n_{0}}|A_{n}^{f_{M}^{\prime }}|>\delta \})<\epsilon $. Applying Lemma \ref%
{lemma_maximal} to $f_{M}^{\prime \prime }=f-f_{M}^{\prime }$ gives $\mu
(\{\sup_{n}|A_{n}^{f_{M}^{\prime \prime }}|>\delta \})<\epsilon $. As a
result, $\mu (\{\sup_{n\geq n_{0}}|A_{n}^{f}|>2\delta \})<2\epsilon $. 
%\end{proof}
$\Box $

\begin{remark}
\label{remarkcompuL1} We remark that a bounded a.e. computable function, as
defined in Definition \ref{a.e.compufunct} is a computable element of $%
L^{1}(X,\mu )$ (see \cite{HoyRojCiE09}). Conversely, if $f$ is a computable
element of $L^{1}(X,\mu )$ then there is a sequence of uniformly computable
functions $f_{n}$ that effectively converge $\mu $-a.e. to $f$.
\end{remark}

\section{Computing pseudorandom points, constructive ergodic theorem\label%
{pseudorandom}}

Let $X$ again be a computable metric space and $\mu $ a computable
probability measure on it. Suppose $X$ is complete.

Points satisfying the above recalled pointwise ergodic theorem for an
observable $f$, will be called\emph{\ typical for }$f$\emph{\ }.

Points which are typical for each $f$ which is continuous with compact
support are called\emph{\ typical for the measure }$\mu $ ( and for the
dynamics).

The set of computable points in $X$ (see Definition \ref{comp_points}) being
countable, is a very small (invariant) set, compared to the whole space. For
this reason, a computable point could be rarely be expected to be typical
for the dynamics, as defined before. More precisely, the Birkhoff ergodic
theorem and other theorems which hold for a full measure set, cannot help to
decide if there exist a computable point which is typical for the dynamics.
Nevertheless computable points are the only points we can use when we
perform a simulation or some explicit computation on a computer.

A number of theoretical questions arise naturally from all these facts. Due
to the importance of forecasting and simulation of a dynamical system's
behavior, these questions also have some practical motivation.

\begin{problem}
Since simulations can only start with computable initial conditions, given
some typical statistical behavior of a dynamical system, is there some
computable initial condition realizing this behavior? how to choose such
points?
\end{problem}

Such points could be called \emph{pseudorandom} points, and a result showing
its existence and their computability from the description of the system
could be seen as a \emph{constructive }version of the pointwise ergodic
theorem.

Meaningful simulations, showing typical behaviors of the dynamics can be
performed if computable, pseudorandom initial conditions exist ( and can be
computed from the description of the system). Thanks to a kind of effective
Borel-Cantelli lemma, in \cite{GHR07} the above problem was solved
affirmatively for a class of systems which satisfies certain technical
assumptions which includes systems whose decay of correlation is faster that 
$C\log ^{2}(time)$. After the results on the estimation of the rate of
convergence given in the previous section we can remove the technical
assumption on the speed of decay of correlations. We will moreover show how
the result is uniform also in $T$ and $\mu $ (the pseudorandom points will
be calculated from a description of the system).

The following result (\cite{GHR07}, Theorem 2 or \cite{GHRP}) shows that
given a sequence $f_{n}$ which converges effectively a.s. to $f$ and given
its speed of convergence then there it is possible to compute points $x_{i}$
for which $f_{n}(x_{i})\rightarrow f(x_{i})$.

\begin{theorem}
\label{theorem_convergence_Borel_Cantelli} Let $X$ be a complete metric
space. Let $f_{n},f$ be uniformly a.e. computable random variables. If $%
f_{n} $ effectively converges almost surely to $f$ then the set $%
\{x:f_{n}(x)\rightarrow f(x)\}$ contains a sequence of uniformly computable
points which is dense in $\mathrm{Supp}(\mu )$. Moreover, the effective
sequence found above depends algorithmically on $f_{n}$ and on the function $%
n(\delta ,\epsilon )$ giving the rate of convergence
\end{theorem}

Since by the results of the previous section, $n(\delta ,\epsilon )$ can be
calculated starting from $T$, $\mu $ and $f$. This result hence can be
directly used to find typical points for the dynamics. Indeed the following
holds (see \cite{GHRP} for the details)

\begin{theorem}
\label{thm1}If $(X,\mu ,T)$ is a computable ergodic system, $f$ is $%
L^{1}(X,\mu )$ and a.e. computable then there is a uniform sequence $x_{i}$
of computable points which is dense on the support of $\mu $ such that for
each $i$ 
\begin{equation*}
\underset{n\rightarrow \infty }{\lim }\frac{1}{n}\sum f(T^{n}(x_{i}))=\int \!%
{f}\,\mathrm{d}{\mu .}
\end{equation*}%
Moreover this sequence can be computed starting from a description of $T$, $%
\mu $ and $f$.
\end{theorem}

Since the above result is uniform in $f$ and it is possible to construct a
r.e. set of computable observables which is dense in the space of compactly
supported continuous functions we can also obtain the following (see again 
\cite{GHRP} for the details)

\begin{theorem}
If $(X,\mu ,T)$ is a computable ergodic system then there is a uniform
sequence $x_{n}$ of computable points which is dense on the support of $\mu $
such that for each $n$, $x_{n}$ is $\mu -$typical. Moreover this sequence
can be computed starting from a description of $T$and $\mu $
\end{theorem}

Now, as an application we get together the uniform results about computation
of invariant measures for the class of piecewise expanding functions shown
in Subsection \ref{PWSEC} to show that in that class the pseudorandom points
can be calculated starting from a description of the map. Indeed since in
that class the physical invariant measure $\mu $ can be calculated starting
from $T$ (see Theorem \ref{pwexpmu}). We hence obtain

\begin{corollary}
Each piecevise expanding map, with computable derivative, as in the
assumptions of Theorem \ref{pwexpmu} has a sequence of pseudorandom points
which is dense in the support of the measure. Moreover this sequence can be
computed starting from the description of the system.
\end{corollary}

Hence establishing some kind of constructive versions of the ergodic theorem.

\section{Conclusions and directions}

In this article we have reviewed some results about the rigorous computation
of invariant measures, invariant sets and typical points. Here, the sentence 
\emph{rigorous computation} means ``computable by a Turing Machine''. Thus,
this can be seen as a theoretical study of which infinite objects in
dynamics can be arbitrarily well approximated by a modern computer (in an
absolute sense), and which cannot. In this line, we presented some general
principles and techniques that allow the computation of the relevant objects
in several different situations. On the other hand, we also presented some
examples in which the computation of the relevant objects is not possible at
all, stating some theoretical limits to the abilities of computers when used
to simulate dynamical systems.

The examples of the second kind, however, seem to be rather rare. An
important question is therefore whether this phenomenon of non-computability
is robust or prevalent in any sense, or if it is rather exceptional. For
example, one could ask whether the non-computability is destroyed by small
changes in the dynamics or whether the non-computability occurs with small
or null probability.

Besides, in this article we have not considered the efficiency (and
therefore the feasibility) of any of the algorithms we have developed. An
important (and more difficult) remaining task is therefore the development
of a \emph{resource bounded} version of the study presented in this paper.
In the case of Julia sets, for instance, it has been shown in \cite{Br1,
Ret, Br2} that hyperbolic Julia sets, as well as some Julia sets with
parabolics, are computable in polynomial time. On the other hand, in \cite%
{BBY06} it was shown that there exists computable Siegel quadratic Julia
sets whose time complexity is arbitrarily high.

For the purpose to compute the invariant measure form the description of the
system, in section \ref{PWSEC} we had to give explicit estimations on the
constants in the Lasota Yorke inequality. This step is important\ also when
techiques different from ours are used (see \cite{L} e.g.). Similar
estimations could be done in other classes of systems, following the way the
Lasota Yorke inequality is proved in each class (although, sometimes this is
not a completely trivial task and requires the \textquotedblleft
effectivization\textquotedblright\ of some step in the proof). A more
general method to have an estimation for the constants or other ways to get
information on the regularity of the invariant measure would be useful.

\end{document}